\documentclass[a4paper,11pt]{amsart}
\usepackage{mathrsfs}
\usepackage{latexsym, amssymb,amsmath,amsthm}
\usepackage[all, knot]{xy}
\xyoption{arc}

\def \dim{\operatorname{dim}}

\def \Hom{\operatorname{Hom}}

\def \C{\mathcal{C}}

\def \Lie{\textrm{Lie}}

\numberwithin{equation}{section}

\newtheorem{theorem}{Theorem}[section]
\newtheorem{lemma}[theorem]{Lemma}
\newtheorem{proposition}[theorem]{Proposition}
\newtheorem{corollary}[theorem]{Corollary}

\newtheorem{example}[theorem]{Example}
\newtheorem{remark}[theorem]{Remark}

\begin{document}
\title[Super cocommutative Hopf algebras]
{Super cocommutative Hopf algebras of finite representation type}
\author{Gongxiang Liu}
\address{Department of Mathematics, Nanjing University, Nanjing 210093, China} \email{gxliu@nju.edu.cn}
\maketitle

\dedicatory{\begin{center}\Small\emph{Dedicated to Professor
Mitsuhiro Takeuchi in honor of his distinguished
career}\end{center}}

\begin{abstract}
Given an algebraically closed field $k$ of characteristic $p>5$, we
classify the finite super algebraic $k$-groups whose algebras of
measures are of finite representation type. Let $\mathscr{G}$ be
such a super group and $\mathscr{\underline{G}}$ the largest
ordinary algebraic $k$-group determined by $\mathscr{G}$. We show
that both $\mathscr{\underline{G}}$ and
$\mathbf{u}(\Lie(\mathscr{G}))$, the restricted enveloping algebra
of Lie superalgebra of $\mathscr{G}$, are of finite representation
type. Moreover, only some special representation-finite algebraic
$k$-groups of dimension zero can appear if $\mathscr{G}\neq
\mathscr{\underline{G}}$. The structure of $\mathscr{G}$ is almost
determined by $\mathscr{\underline{G}}$ and
$\mathbf{u}(\Lie(\mathscr{G}))$. The Auslander-Reiten quivers are
determined by showing they are Nakayama algebras. \vskip 5pt

\noindent{\bf Keywords} \ \ Super group scheme, Representation type, Nakayama algebra  \\
\noindent{\bf 2000 MR Subject Classification} \ \ 14L15, 16G60,
16W30, 17B50
\end{abstract}

\section{Introduction}
According to the fundamental result of Drozd  \cite{Dr}, every
finite dimensional algebra exactly belongs to one of following three
kinds of algebras: algebras of finite representation type, algebras
of tame type and wild algebras. For the algebras of the former two
kinds, a classification of indecomposable modules seems feasible. By
contrast, the module category of a wild algebra, being
``complicated" at least as that of any other algebra, can't afford
such a classification. Inspired by the Drozd's result, one is often
interested in classifying a given kind of algebras according to
their representation type.

This paper is concerned with the classification of
(representation-finite) super cocommutative Hopf algebras over
algebraically closed fields of positive characteristic. It is known
that such an algebra can be viewed as the group algebras of a finite
super algebraic $k$-group. Special cases are group algebras
associated to finite algebraic $k$-groups, that is,
finite-dimensional cocommutative Hopf algebras, as well as
restricted enveloping algebras of restricted Lie superalgebras. The
representation theory of both of these classes has received
considerable attention. The very detailed information on the
structure of representation-finite and tame cocommutative Hopf
algebras, partially considered by the pioneers as Hochschild
\cite{Hoc}, Feldvoss and Strade \cite{FS}, Pfautsch and Voigt
\cite{PV,Voi}, etc, has been ultimately gotten by Farnsteiner and
his corporators continuous \cite{Far2,Far3,Far4,FV}. Also, the
restricted Lie superalgebras of finite representation type were
classified by Farnsteiner too \cite{Far1}. Our final goal will be
the extension of these results to arbitrary super cocommutative Hopf
algebras.

There are two ways to connect super Hopf algebras $H$ with ordinary
Hopf algebras and both of them will be used freely in the paper. One
is the Radford-Majid bosonization \cite{Majid,Radford}, which
constructs from $H$ an ordinary Hopf algebra $H\rtimes
k\mathbb{Z}_{2}$. Another one, given by Masuoka \cite{Mas}, states
that if $H$ is super cocommutative, there is a unit-preserving
isomorphism
$$H\cong \underline{H}\otimes \wedge(V_{H})$$
as super left $\underline{H}$-module coalgebras, where
$\underline{H}$ is the largest ordinary sub Hopf algebra and
$V_{H}=P(H)_{1}$. These two ways will be recalled in Section 2.

So, in philosophy, one just need to know ``how to" reduce the
research of  representation theory of a super cocommutative Hopf
algebras $H$ to that of $\underline{H}$ and $V_{H}$.  Both Section 3
and Section 4 are designed to give methods of such reduction. The
main result of Section 3 shows that $\textbf{u}(P(H))$, which
controls $V_{H}$ essentially, has finite representation type
provided that $H$ is of finite representation type. Under assumption
that $H$ is of finite representation type, the structure of
$\underline{H}$ are shown to be quite special. We will see in
Section 4 that either $\underline{H}$ is semisimple or the
$V$-uniserial group attached to it has height $\leq 1$. Due to the
lack of Mackey decomposition for super algebraic groups, one has to
apply other methods. It turns out the concept of complexity, which
are shown effective in dealing with infinitesimal groups, is also
quite useful in our case. And in section 2, some notions and
computations relevant to our purpose, particularly the concept of a
path coalgebra and complexity,
 are summarized. Combining the results gotten
in Sections 3,4, the representation-finite super finite algebraic
groups are determined in Section 5. The representation theory of
them are determined through showing they are always Nakayama
algebras in the last section.

\section{preliminaries}
Throughout we will be working over a field $k$. All spaces are
$k$-spaces. For short, $\otimes_{k}$ is just denoted by $\otimes$.

\subsection{Path coalgebras.}
Given a quiver
 $Q=(Q_{0},Q_{1})$ with $Q_{0}$ the set of vertices and $Q_{1}$
 the set of arrows, denote by $kQ$ and $kQ^{c}$, the
 $k$-space with basis the set of all paths in $Q$ and the path coalgebra of $Q$, respectively. Note
 that they are all graded with respect to length grading. For
 $\alpha \in Q_{1}$, let $s(\alpha)$ and $t(\alpha)$ denote
 respectively the starting and ending vertex of $\alpha$.

 Recall that the comultiplication of the path coalgebra $kQ^{c}$
 is defined by
 $$\Delta(p)=\alpha_{l}\cdots \alpha_{1}
 \otimes s(\alpha_{1})+\sum_{i=1}^{l-1}\alpha_{l}\cdots \alpha_{i+1}
 \otimes\alpha_{i}\cdots \alpha_{1}+ t(\alpha_{l})\otimes \alpha_{l}\cdots \alpha_{1}$$
 for each path $p=\alpha_{l}\cdots \alpha_{1}$ with each $\alpha_{i}\in
 Q_{1}$; and $\varepsilon(p)=0$ for $l\geq 1$ and $1$ if $l=0$
($l=0$ means $p$ is a vertex). This is a pointed coalgebra.

For a quiver $Q$, define $$kQ_{d}:=\oplus_{i=0}^{d-1}kQ(i)$$ where
$Q(i)$ is the set of all paths of length $i$ in $Q$. Our interested
quiver is the simplest one, a loop $\circlearrowleft$. For any
natural number $n$, denote the unique path of length $n$ of
$k\circlearrowleft$ by $\alpha_{n}$. In particular,
$k\circlearrowleft_{p^{n}}$ has a basis
$1,\alpha_{1},\alpha_{2},\ldots,\alpha_{p^{n}-1}$.

\subsection{Representation type.} A finite-dimensional algebra $A$ is said to be of
\emph{finite representation type} provided there are finitely many
non-isomorphic indecomposable $A$-modules. $A$ is of \emph{tame
type} or $A$ is a \emph{tame} algebra if $A$ is not of finite
representation type, whereas for any dimension $d>0$, there are
finite number of $A$-$k[T]$-bimodules $M_{i}$ which are free of
finite rank as right $k[T]$-modules such that all but a finite
number of indecomposable $A$-modules of dimension $d$ are isomorphic
to $M_{i}\otimes_{k[T]}k[T]/(T-\lambda)$ for $\lambda\in k$. We say
that $A$ is of \emph{wild type} or $A$ is a \emph{wild} algebra  if
there is a finitely generated $A$-$k\langle X,Y\rangle $-bimodule
$B$ which is free as a right $k\langle X,Y\rangle $-module such that
the functor $B\otimes_{k\langle X,Y\rangle}-\;\;$ from mod-$k\langle
X,Y\rangle$, the category of finitely generated $k\langle
X,Y\rangle$-modules, to mod-$A$, the category of finitely generated
$A$-modules, preserves indecomposability and reflects isomorphisms.
See \cite{E} for more details.

\subsection{Super cocommutative Hopf algebras.} We recall the two ways
connecting super cocommutative Hopf algebras with usual Hopf
algebras in this subsection. Let $J$ be a Hopf algebra with
bijective antipode and $^{J}_{J}\mathscr{YD}$ the category of the
Yetter-Drinfeld modules with left $J$-module action and left
$J$-comodule coaction. It is naturally forms a braided monoidal
category with the braiding
$$c_{M,N}:\;M\otimes N\to N\otimes M,\;\;m\otimes n\mapsto \sum n_{0}\otimes n_{-1}\cdot m,$$
where $n\mapsto \sum n_{-1}\otimes n_{0},\;N\to J\otimes N$ denotes
the comodule structure, as usual. Let $A$ be a Hopf algebra in
$^{J}_{J}\mathscr{YD}$. In particular, $A$ is a left $J$-module
algebra and left $J$-comodule coalgebra. The Radford-Majid
bosonization \cite{Majid,Radford} gives rise to an ordinary Hopf
algebra, $A\rtimes J$. As an algebra, this is the smash product $A\#
J$, and it is the smash coproduct as a coalgebra. In particular, a
super Hopf algebra $H$ is a Hopf algebra in
$^{k\mathbb{Z}_{2}}_{k\mathbb{Z}_{2}}\mathscr{YD}$ (see Section 2 in
\cite{Mas}) and hence we get a usual Hopf algebra $H\rtimes
k\mathbb{Z}_{2}$. The following result is a direct consequence of
Proposition 3.2 in Chapter VI of \cite{ARS}.

\begin{lemma} Assume that characteristic of $k$ is not $2$. Then $H$ and $H\rtimes
k\mathbb{Z}_{2}$ have the same representation type.
\end{lemma}

For a super Hopf algebra $H=H_{0}\oplus H_{1}$, apart from its
ordinary representations, one also can consider its super
representations. That is, the $\mathbb{Z}_{2}$-graded $H$-modules.
Clearly, super representations of $H$ are just the ordinary
representations of the ordinary Hopf algebra $H\rtimes
k\mathbb{Z}_{2}$. Thus,

\begin{lemma} Let $H$ be a super Hopf algebra. Then the category of super
$H$-modules is equivalent to the category of $H\rtimes
k\mathbb{Z}_{2}$-modules.
\end{lemma}

Combining Lemma 2.1 and Lemma 2.2, if  characteristic of $k$ is not
2, the representation type of $H$ as an ordinary algebra is indeed
the same with that of $H$ when we consider it as a super algebra. In
this paper, we will always consider the ordinary representations
except in the proof of Theorem 6.1.

An algebra $A$ is a \emph{Nakayama algebra }if each indecomposable
$A$-module is uniserial. The following lemma is the Theorem 2.14 in
Chapter IV of \cite{ARS}.

\begin{lemma} Let $G$ be a finite group such that $|G|$ is
invertible in $k$ and $A$ is a finite-dimensional $kG$-module
algebra. Then $A\# kG$ is a Nakayama algebra if and only if $A$ is
so.
\end{lemma}
Let $H$ be a super Hopf algebra, we call $H$ a \emph{super Nakayama
algebra} if each super indecomposable $H$-module is uniserial. Owing
to Lemma 2.2 and Lemma 2.3, $H$ is Nakayama if and only if it is
super Nakayama. This fact will be used in the proof of Theorem 6.1.

Now let $H=H_{0}\oplus H_{1}$ be a super cocommutative Hopf algebra
over $k$. Define
$$\underline{H}:=\Delta^{-1}(H_{0}\otimes H_{0}).$$
This is the largest ordinary cocommutative sub Hopf algebra of $H$.
Denote the set of primitives in $H$ by $P(H)$ and define
$$V_{H}:=P(H)_{1}$$
the vector space of odd primitives in $H$. Choose a totally ordered
$k$-basis $X=(x_{\lambda})_{\lambda}$ of $V_{H}$. Then,
$x_{\lambda}\wedge x_{\mu}\wedge\cdots\wedge x_{\nu}
(x_{\lambda}<x_{\mu}<\cdots<x_{\nu})$ form a $k$-basis of $\wedge
(V_{H})$, and $x_{\lambda}\wedge x_{\mu}\wedge\cdots\wedge
x_{\nu}\mapsto x_{\lambda} x_{\mu}\cdots x_{\nu}$ gives a
unit-preserving super coalgebra map from $\wedge (V_{H})$ to $H$. We
collects some facts about $H$, which were given essentially by
Masuoka in \cite{Mas}, as follows.
\begin{lemma} \emph{(1)} The induced left
$\underline{H}$-linear map
$$\phi=\phi_{X}:\;\underline{H}\otimes \wedge
(V_{H})\longrightarrow H$$ is a unit-preserving isomorphism of super
left $\underline{H}$-module coalgebra.

\emph{(2)} As an algebra, $H$ is generated by $\underline{H}$ and
$V_{H}$.

\emph{(3)} $V_{H}$ is a right $\underline{H}$-module under the
conjugation $v\cdot h:=\sum S(h_{(1)})vh_{(2)}$ and
$$vh=\sum h_{(1)}(v\cdot h_{(2)})$$
for $v\in V_{H}$ and $h\in \underline{H}$.

\emph{(4)} For any $u,v\in V_{H}$, we have $uv+vu\in P(H)$.
\end{lemma}
\begin{proof} (1) and the first part of (3) are Theorem 3.6 and Proposition
3.9 (1) in \cite{Mas} respectively. (2) is a direct consequence of
Proposition 3.9 (2) in \cite{Mas}. Both (4) and the second part of
(3) can be gotten easily by direct computations.
\end{proof}

\noindent{\textbf{Convention.}} Due to (1) and (2) of above lemma,
sometimes we use the notation $\underline{H}\langle V_{H}\rangle$ to
denote the super cocommutative Hopf algebra $H$. This is convenient.
For example, let $K\subset \underline{H}$ be a sub Hopf algebra
containing $P(H)$ and $K'$ the sub super cocommutative Hopf algebra
generated by $K$ and some $V\subset V_{H}$. Then we have
$\underline{K'}=K$ (by (4) of above lemma) and $V_{K'}=V$. So
$K'=K\langle V\rangle$. Moreover, if $\dim_{k}V_{H}=1$, then we will
simply use the notion $\underline{H}\langle v\rangle$ instead of
$\underline{H}\langle V_{H}\rangle$
for any nonzero element $v\in V_{H}$. \\

Let $C$ be a (super) coalgebra, define
$C^{+}:=\textrm{Ker}(\varepsilon)$ as usual.
\begin{lemma} Let $K\subset \underline{H}$ be a sub normal Hopf algebra
containing $P(H)$ and $V\subset V_{H}$ a subspace of $V_{H}$. Then
there is a Hopf isomorphism
$$\underline{H}/K^{+}\underline{H} \cong \underline{H}\langle V\rangle/(K\langle V\rangle)^{+}\underline{H}\langle V\rangle.$$
\end{lemma}
\begin{proof} By $K$ is normal, we have exact sequence of Hopf
algebras $K\hookrightarrow \underline{H}\twoheadrightarrow
\underline{H}/K^{+}\underline{H}$. Also, we have an obvious exact
sequence $V\hookrightarrow V\twoheadrightarrow 0$. Owing to Theorem
3.13 (3) in \cite{Mas}, the sequence $K\langle V\rangle
\hookrightarrow \underline{H}\langle V\rangle\twoheadrightarrow
\underline{H}/K^{+}\underline{H}$ is also exact. Thus the conclusion
is proved.
\end{proof}

\subsection{Complexity.} Let $A$ be an associative
algebra, $M$ an $A$-module with minimal projective resolution
$$\cdots \rightarrow P_{n}\rightarrow P_{n-1}\rightarrow \cdots \rightarrow P_{0}\rightarrow M\rightarrow 0.$$
Then the \emph{complexity} of $M$ is defined to be the integer
$$\C_{A}(M):=\textrm{min}\{c\in \mathbb{N}_{0}\cup \infty\;| \exists \lambda>0:\textrm{dim}_{k}P_{n}\leq
\lambda n^{c-1},\;\forall \;n\geq 1\}.$$

For our purpose, we need consider the following examples.

\begin{example}
\emph{(1) Let $A$ be a self-injective algebra of finite
representation type, then it is well-known that $\C_{A}(M)\leq 1$
for any $A$-modules $M$. }

\emph{(2) Consider the algebra $A=k[x,y]/(x^{n},y^{2})$ for some
$n>1$. It is a local algebra and we denote the unique simple module
by $k$. We can construct the minimal projective resolution of $k$ as
follows.}

\begin{figure}[hbt]
\begin{picture}(100,80)(-170,-30)
\put(0,0){\makebox(0,0){$ k$}} \put(-30,0){\vector(1,0){22}}
\put(5,0){\vector(1,0){10}} \put(22,0){\makebox(0,0){$ 0$}}

\put(-40,0){\makebox(0,0){$A$}} \put(-80,0){\vector(1,0){30}}
 \put(-100,0){\makebox(0,0){$A\oplus
A$}} \put(-70,30){\makebox(0,0){$\Omega (k)$}}
\put(-65,20){\vector(1,-1){15}} \put(-90,5){\vector(1,1){15}}

\put(-150,0){\vector(1,0){30}}
 \put(-180,0){\makebox(0,0){$A\oplus A\oplus A$}}
\put(-140,-30){\makebox(0,0){$\Omega^{2} (k)$}}
\put(-125,-25){\vector(1,1){15}}\put(-175,-10){\vector(1,-1){15}}

 \put(-245,0){\vector(1,0){30}} \put(-260,0){\makebox(0,0){$\cdots$}}
 \put(-230,30){\makebox(0,0){$\Omega^{3} (k)$}}
 \put(-225,20){\vector(1,-1){15}} \put(-250,5){\vector(1,1){15}}

\end{picture}
\end{figure}

\emph{Here $\Omega(M)$ is the kernel of a minimal projective cover
of the $A$-module $M$. It is not hard to show that}

\begin{figure}[hbt]
\begin{picture}(100,100)(-20,20)

\put(-50,80){\makebox(0,0){$\Omega^{2}(k)=$}}
\put(50,115){\makebox(0,0){$(y,-x)$}}
\put(50,110){\makebox(0,0){$\bullet$}}
\put(50,110){\line(-1,-1){20}}
\put(40,100){\makebox(0,0){$\bullet$}}
\put(30,90){\makebox(0,0){$\bullet$}}
\put(25,85){\makebox(0,0){$\cdot$}}
\put(20,80){\makebox(0,0){$\cdot$}}
\put(15,75){\makebox(0,0){$\cdot$}}

\put(10,70){\makebox(0,0){$\bullet$}}

 \put(10,70){\line(-1,-1){20}}
\put(0,60){\makebox(0,0){$\bullet$}}
\put(-10,50){\makebox(0,0){$\bullet$}}
\put(-10,40){\makebox(0,0){$x^{n-1}y$}}

\put(-20,60){\line(1,-1){10}}\put(-20,60){\makebox(0,0){$\bullet$}}\put(-20,65){\makebox(0,0){$x^{n-1}$}}

\put(50,110){\line(1,-1){10}} \put(60,100){\makebox(0,0){$\bullet$}}
\put(40,100){\line(1,-1){10}}

\put(75,115){\makebox(0,0){$y$}}
\put(70,110){\makebox(0,0){$\bullet$}}
\put(70,110){\line(-1,-1){20}}
\put(60,100){\makebox(0,0){$\bullet$}}\put(50,90){\makebox(0,0){$\bullet$}}
\put(40,80){\makebox(0,0){$\cdot$}}
\put(35,75){\makebox(0,0){$\cdot$}}
\put(30,70){\makebox(0,0){$\cdot$}} \put(20,60){\line(-1,-1){10}}
\put(20,60){\makebox(0,0){$\bullet$}}
\put(10,50){\makebox(0,0){$\bullet$}}

\put(10,70){\line(1,-1){10}} \put(0,60){\line(1,-1){10}}
\put(27,50){\makebox(0,0){$x^{n-1}y$}}
\end{picture}
\end{figure}
\emph{Here we use the line $/$ to denote the action of $x$ and
$\backslash$ the action of $y$. Define}

\begin{figure}[hbt]
\begin{picture}(80,80)(-20,30)
\put(-50,80){\makebox(0,0){$N:=$}} \put(50,115){\makebox(0,0){$x$}}
\put(50,110){\makebox(0,0){$\bullet$}}
\put(50,110){\line(-1,-1){20}}
\put(40,100){\makebox(0,0){$\bullet$}}
\put(30,90){\makebox(0,0){$\bullet$}}
\put(25,85){\makebox(0,0){$\cdot$}}
\put(20,80){\makebox(0,0){$\cdot$}}
\put(15,75){\makebox(0,0){$\cdot$}}

\put(10,70){\makebox(0,0){$\bullet$}}

 \put(10,70){\line(-1,-1){10}}
\put(0,60){\makebox(0,0){$\bullet$}}
\put(-5,55){\makebox(0,0){$x^{n-1}$}}

\put(50,110){\line(1,-1){10}} \put(60,100){\makebox(0,0){$\bullet$}}
\put(40,100){\line(1,-1){10}}

\put(65,105){\makebox(0,0){$xy$}} \put(60,100){\line(-1,-1){10}}
\put(60,100){\makebox(0,0){$\bullet$}}\put(50,90){\makebox(0,0){$\bullet$}}
\put(40,80){\makebox(0,0){$\cdot$}}
\put(35,75){\makebox(0,0){$\cdot$}}
\put(30,70){\makebox(0,0){$\cdot$}} \put(20,60){\line(-1,-1){10}}
\put(20,60){\makebox(0,0){$\bullet$}}
\put(10,50){\makebox(0,0){$\bullet$}}

\put(10,70){\line(1,-1){10}} \put(0,60){\line(1,-1){10}}
\put(27,50){\makebox(0,0){$x^{n-1}y$}}

\end{picture}
\end{figure}

\emph{Then through direct computations, we have for any $3\leq i\in
\mathbb{N}$}

$$\Omega^{i}(k)=\left \{
\begin{array}{ll}  \Omega^{i-1}(k)\oplus N& \;\;\;\;\emph{if}\;i\;\emph{is\;odd,}\\
\Omega^{i-1}(k)\oplus \Omega(N) &
\;\;\;\;\emph{if}\;i\;\emph{is\;even.}
\end{array}\right. $$

\emph{By this, we indeed have get
$$P_{n}\cong A^{(n+1)}$$
for a minimal projective resolution $P_{\bullet}\to k$. This implies
that $\C_{A}(k)=2$.}

\end{example}

\section{Structure of  $\mathbf{u}(\Lie(\mathscr{G}))$}

Throughout this and following sections, we assume that $k$ is an
algebraically closed field with characteristic $p>5$. Let
$\mathscr{G}$ be a finite super algebraic $k$-group and
$H(\mathscr{G})$ be its algebra of measures. That is,
$H(\mathscr{G})=(\mathcal{O}(\mathscr{G}))^{\ast}$. Then
$H(\mathscr{G})$ is a finite-dimensional super cocommutative Hopf
algebra and $\Lie(\mathscr{G})=P(H(\mathscr{G}))$ is a restricted
Lie superalgebra. Denote $\mathbf{u}(\Lie(\mathscr{G}))$ the
restricted enveloping algebra of $\Lie(\mathscr{G})$ and it is a sub
super Hopf algebra of $H(\mathscr{G})$. The purpose of this section
is to show that $\mathbf{u}(\Lie(\mathscr{G}))$ is of finite
representation type provided $H(\mathscr{G})$ is so.

Let $(L_{0},[p])$ be a restricted Lie algebra.  An element $x\in
L_{0}$ is called to be \emph{toral} if $x^{[p]}=x$ and
\emph{$p$-nilpotent} if there exists some $n\in \mathbb{N}$ such
that $x^{[p]^{n}}=0$. We denote by $N(L_{0})$, $T(L_{0})$ and
$C(L_{0})$, the largest nilpotent ideal, toral ideal and center of
$L_{0}$. If $X\subset L_{0}$ is a subset, then $X_{p}$ denotes the
$p$-subalgebra of $L_{0}$ that is generated by $X$. See \cite{SF}
for details. The following conclusion was given in \cite{Far1} as
Lemma 4.1 and Theorem 4.2.

\begin{lemma} \emph{(1)} Let $V$ be an $k$-vector space with exterior
algebra $\wedge(V)$. If $\mathcal{C}_{\wedge(V)}(k)\leq 1$, then
$\dim_{k}V\leq 1$.

\emph{(2)} Let $L=L_{0}\oplus L_{1}$ be a restricted Lie
superalgebra with $L_{1}\neq 0$. Then $\emph{\textbf{u}}(L)$ has
finite representation type if and only if there exists a toral
element $t_{0}\in L_{0}$, a $p$-nilpotent element $x_{0}\in L_{0}$,
and $y\in L_{1}$ such that $L=L_{0}\oplus ky$, $(kx_{0})_{p}\subset
[L_{1},L_{1}]_{p}\subset N(L_{0}),
\;L_{0}=N(L_{0})+kt_{0},\;N(L_{0})=T(L_{0})\oplus (kx_{0})_{p}$.
\end{lemma}

\emph{In the following of this paper, we fix the notion $x_{0}$ and
$y$ to denote the elements given in this lemma.} We can prove our
conclusion now.

\begin{proposition} Let $\mathscr{G}$ be a finite super
$k$-group and $H(\mathscr{G})$ be its algebra of measures. If
$H(\mathscr{G})$ is of finite representation type then
$\mathbf{u}(\emph{Lie}(\mathscr{G}))$ is so too.\end{proposition}

The $(1)\Rightarrow (2)$ part of the proof of Theorem 4.2 in
\cite{Far1} can be essentially applied to our case except some
delicate points. So our proof looks like ``cut and paste". For
safety and convenience of readers, we still write it out.

\begin{proof} For simplicity, denote the restricted Lie superalgebra $\Lie(\mathscr{G})$ by
$L=L_{0}+L_{1}$. If $L_{1}=0$, the conclusion can be proved easily.
In fact, Example 2.6 (1) implies that
$\mathcal{C}_{H(\mathscr{G})}(k)\leq 1$ where we consider $k$ as the
trivial $H(\mathscr{G})$ module through the map $\varepsilon:\;
H(\mathscr{G})\to k$. Since $H(\mathscr{G})$ is a free
$\textbf{u}(L)$-module, $\mathcal{C}_{\textbf{u}(L)}(k)\leq 1$.
Consequently, $\textbf{u}(L)$ has finite representation type by
Theorem 2.4 in \cite{FS}. Therefore, one can assume that $L_{1}\neq
0$. Our goal is to show that $L$ indeed has the structure as
described in Lemma 3.1 (2). We divide the task into several steps.

(a) \emph{{There exists a toral element  $t_{0}\in L_{0}$, a
$p$-nilpotent element $x_{0}\in L_{0}$ such that
$L_{0}=N(L_{0})+kt_{0},\;N(L_{0})=T(L_{0})\oplus (kx_{0})_{p}$.}}
Owing to the discussion in the above paragraph, $\textbf{u}(L_{0})$
has finite representation type and thus Theorem 4.3 in \cite{Far2}
implies the desired result.

(b) \emph{Let $T:=T(L_{0})+kt_{0}$. Then $T$ is a maximal torus of
$L_{0}$ and there exists at most one root $\alpha$ relative to $T$.
The corresponding root space $(L_{0})_{\alpha}$ has dimension 1.}
Totally the same with the part $(1)\Rightarrow (2)$ (c) of the proof
of Theorem 4.2 in \cite{Far1}.

In the following of the proof, we decompose the $T$-module $L_{1}$
into its weight spaces and write $L_{1}=\bigoplus_{\lambda\in
\mathscr{W}}(L_{1})_{\lambda}$, where $\mathscr{W}\subset T^{\ast}$
is the set of weights of $L_{1}$ relative to $T$.

(c) \emph{Let $\beta\in \mathscr{W} \backslash
\{0,\frac{1}{2}\alpha\}$. Then $\dim_{k}(L_{1})_{\beta}=1$,
$[(L_{1})_{\beta},(L_{1})_{\beta}]=0$, and $\mathscr{W}\subset
\{0,\frac{1}{2}\alpha,\beta,-\beta\}$ or $\mathscr{W}\subset
\{0,\frac{1}{2}\alpha,\beta,\alpha-\beta\}$.} See the part
$(1)\Rightarrow (2)$ (d) of the proof of Theorem 4.2 in \cite{Far1}.

(d) \emph{Suppose that $\mathscr{W} \backslash
\{0,\frac{1}{2}\alpha\}\neq \emptyset$. Then there exists $\gamma\in
\mathscr{W} \backslash \{0,\frac{1}{2}\alpha\}$ such that
$[(L_{0})_{\alpha}, (L_{1})_{\gamma}]=0$.} See the part
$(1)\Rightarrow (2)$ (e) of the proof of Theorem 4.2 in \cite{Far1}.

(e) \emph{If $\mathscr{W} \backslash \{0,\frac{1}{2}\alpha\}\neq
\emptyset$, then $L=T(L_{0})\oplus ky_{1}$, where $ky_{1}=L_{1}$,
and $[y_{1},y_{1}]=0$.}

By (d), there exists $\gamma\in \mathscr{W} \backslash
\{0,\frac{1}{2}\alpha\}$ such that $[(L_{0})_{\alpha},
(L_{1})_{\gamma}]=0$. According to (b),
$(L_{0})_{\alpha}=kx_{\alpha}$, for some $x_{\alpha}\in L_{0}$, such
that $N(L_{0})=C(L_{0})+kx_{\alpha}$. Since
$N(L_{0})=C(L_{0})+kx_{0}$ (by the structure of $L_{0}$ described in
(a)), it follows that $x_{0}=x_{a}+z$ for $z\in C(L_{0})\subset
(L_{0})_{0}$. This implies
$[x_{0},(L_{1})_{\gamma}]=[z,(L_{1})_{\gamma}]\subset
(L_{1})_{\gamma}$. As $x_{0}$ is $p$-nilpotent and
$\dim_{k}(L_{1})_{\gamma}=1$, $[x_{0},(L_{1})_{\gamma}]=0$. Consider
the restricted Lie superalgebra $\mathscr{L}:=(kx_{0})_{p}\oplus
(L_{1})_{\gamma}$, then its restricted enveloping algebra is
$\textbf{u}(k(x_{0})_{p})\otimes \wedge ((L_{1})_{\gamma})$. Also,
since $H(\mathscr{G})$ is projective over $\textbf{u}(\mathscr{L})$,
$\mathcal{C}_{\textbf{u}(\mathscr{L})}(k)\leq 1$. Then the
K$\ddot{\textrm{u}}$nneth formula implies
$\textbf{u}(\mathscr{L})=k$ and so $x_{0}=0$. Therefore,
$N(L_{0})=T(L_{0})=L_{0}$.

From (c) we now obtain $\mathscr{W}\subset \{0,\gamma,-\gamma\}$ and
$L=T(L_{0})\oplus (L_{1})_{0}\oplus (L_{1})_{\gamma}\oplus
(L_{1})_{-\gamma}$. Let $X:=Ker\gamma$, it is a $p$-ideal of
$L_{0}$. Applying the same computation used in the part
$(1)\Rightarrow (2)$ (f) of the proof of Theorem 4.2 in \cite{Far1},
$X$ satisfies the relations $[X,L_{1}]=0$ and $[L_{1},L_{1}]\subset
X$. Let $\mathscr{L}:=L/X$. Consider the quotient super Hopf algebra
$H(\mathscr{G})/(X)$, and of course it is of finite representation
type. Owing to
 $H(\mathscr{G})$ is faithfully flat over $\textbf{u}(\mathscr{L})$
 (in fact, $H(\mathscr{G})$ is free over $\textbf{u}(\mathscr{L})$),
 $\textbf{u}(\mathscr{L})\cap (X)=X$. Then $\textbf{u}(\mathscr{L})$ is a sub
 super Hopf algebra of $H(\mathscr{G})/(X)$ and so $\mathcal{C}_{\textbf{u}(\mathscr{L})}(k)\leq
 1$. Therefore, we also have $\mathcal{C}_{\wedge \mathscr{L}_{1}}(k)\leq
 1$, and Lemma 3.1 (1) yields $\dim_{k}L_{1}=\dim_{k}\mathscr{L}_{1}\leq
 1$.

 In view the result in (e), we shall hence forth assume the
 $\frac{1}{2}\alpha$ and 0 are the only weights of $L_{1}$ relative
 $T$.

 (f) \emph{$[N(L_{0}),L_{1}]=0$.} See the part
$(1)\Rightarrow (2)$ (g) of the proof of Theorem 4.2 in \cite{Far1}.

(g) \emph{$[L_{1},L_{1}]_{p}\subset N(L_{0})$ and
$\dim_{k}L_{1}=1$.}

It follows from (f) that $[N(L_{0}),[L_{1},L_{1}]]=0$, proving that
$[L_{1},L_{1}]$ is contained in the centralizer of
$C_{L_{0}}(N(L_{0}))$ which equals to $N(L_{0})$ by (a). Thus
$[L_{1},L_{1}]_{p}\subset N(L_{0})$. Let $\mathscr{L}:=L/N(L_{0})$.
Also, consider the quotient $H(\mathscr{G})/(N(L_{0}))$ and it has
finite representation type. Similar to the proof of part (e),
$\textbf{u}(\mathscr{L})$ is a sub
 super Hopf algebra of $H(\mathscr{G})/(N(L_{0}))$ and so $\mathcal{C}_{\textbf{u}(\mathscr{L})}(k)\leq
 1$. Therefore, we also have $\mathcal{C}_{\wedge \mathscr{L}_{1}}(k)\leq
 1$, and Lemma 3.1 (1) yields $1\leq \dim_{k}L_{1}=\dim_{k}\mathscr{L}_{1}\leq
 1$.

 (h) $(kx_{0})_{p}\subset [L_{1},L_{1}]_{p}$.

 By (g), $L_{1}=ky$ for some $0\neq y\in L_{1}$. Put $v:=[y,y]$.
 Owing to (f) and (g), $[L_{1},L_{1}]_{p}=(kv)_{p}$ is an ideal of
 $L$. Let $\mathscr{L}:=L/(kv)_{p}$ and consider the quotient
 $H(\mathscr{G})/((kv)_{p})$. Using the methods developed in proof
 of  (e) again, one has $\mathcal{C}_{\textbf{u}(\mathscr{L})}(k)\leq
 1$. Denote the natural projection $L\to \mathscr{L}$ by $\pi$.
 Since $[\mathscr{L}_{1},\mathscr{L}_{1}]=0$ and $(kx_{0})_{p}$
 operates trivially on $L_{1}$, it follows that $\textbf{u}((k\pi(x_{0}))_{p}\oplus
 \mathscr{L}_{1})$ is isomorphic to $\textbf{u}((k\pi(x_{0}))_{p})\otimes
 k[X]/(X^{2})$. By $\mathcal{C}_{\textbf{u}(\mathscr{L})}(k)\leq
 1$ and $\textbf{u}(\mathscr{L})$ is projective over  $\textbf{u}((k\pi(x_{0}))_{p}\oplus
 \mathscr{L}_{1})$, $\mathcal{C}_{\textbf{u}((k\pi(x_{0}))_{p}\oplus
 \mathscr{L}_{1})}(k)\leq 1$. The K$\ddot{\textrm{u}}$nneth formula
 implies that $\textbf{u}((k\pi(x_{0}))_{p})=k$ and consequently $(kx_{0})_{p}\subset
 [L_{1},L_{1}]_{p}$.

 By the results getting in (a)-(h), $L$ has the structure
 described in sufficiency's part of Lemma 3.1 (2). Thus
 $\textbf{u}(L)$ has finite representation type.
 \end{proof}

\section{Structure of $\underline{H}(\mathscr{G})$}

 Recall in Section 3, for a super
cocommutative Hopf algebra $H$ we denote its largest cocommutative
sub Hopf algebra by $\underline{H}$. Let $\mathscr{G}$ be a finite
super algebraic $k$-group and $H(\mathscr{G})$ be its algebra of
measures. Denote $\mathscr{\underline{G}}$ the largest ordinary
algebraic $k$-group of $\mathscr{G}$, i.e., by definition its
algebra of measures $H(\mathscr{\underline{G}})$ is
$\underline{H}(\mathscr{G})$. That is,
$H(\mathscr{\underline{G}})=\underline{H}(\mathscr{G})$.
\emph{Throughout this section, we always assume that
$\mathscr{G}\neq \mathscr{\underline{G}}$.} The task of this section
is to analysis the structure of $\mathscr{\underline{G}}$.

\begin{proposition} Assume that $H(\mathscr{G})$ has finite
representation type. Then $\underline{H}(\mathscr{G})$ has finite
representation type too.
\end{proposition}
\begin{proof} Denote $H(\mathscr{G})$ by $H$ for simplicity.  Owing to Proposition 3.2 and Lemma 3.1 (2),
$\dim_{k}\Lie(\mathscr{G})_{1}=1$. So there exists $0\neq y\in
V_{H}$ (Recall in Section 2, $V_{H}$ was defined to be $P(H)_{1}$)
such that $V_{H}=ky$. By Lemma 2.4 (3), $V_{H}$ is right
$\underline{H}$-module. Thus, there exists an algebra map
$$\chi:\;\underline{H}\longrightarrow k$$
such that $y\cdot h=\chi(h)y$ for $h\in \underline{H}$. Let
$\alpha\in \textrm{Aut}(\underline{H})$ be the algebra automorphism
determined by $\chi$, that is, $\alpha(h):=(id\ast \chi)(h)=\sum
h_{(1)}\chi(h_{(2)})$ for $h\in \underline{H}$. By Lemma 2.4 (3), we
always have for $h\in \underline{H}$
\begin{equation} yh=\alpha(h)y.
\end{equation}
Thus $$H=\underline{H}\oplus \underline{H}y$$ as
$\underline{H}$-bimodules and Lemma 3.1 (a) in Chapter VI of
\cite{ARS} implies that $\underline{H}$ is of finite representation
type too.
\end{proof}

In the following of this section, we always assume that
$H(\mathscr{G})$ has finite representation type. By the proof of
this proposition there exists $0\neq y\in V_{H}$ such that
$H=\underline{H}\langle y\rangle$ (See the Convention after Lemma
2.4). It is known any ordinary finite algebraic $k$-group
$\mathscr{H}$ can be decomposed into a semidirect product
$\mathscr{H}=\mathscr{H}^{\circ}\rtimes \mathscr{H}_{red}$ with a
constant group $\mathscr{H}_{red}$ and a normal infinitesimal
subgroup $\mathscr{H}^{\circ}$. In particular,
$$\mathscr{\underline{G}}=\mathscr{\underline{G}}^{\circ}\rtimes \mathscr{\underline{G}}_{red}.$$
 With such notions,

 \begin{lemma} $H(\mathscr{\underline{G}}_{red})$ is always
 semisimple.
 \end{lemma}
 \begin{proof} At first, assume that $[y,y]=0$. If $H(\mathscr{\underline{G}}_{red})$ is not semisimple, then there
exists $g\in \mathscr{\underline{G}}_{red}$ of order $p$. Since the
automorphism group of $ky$ is the multiplicative group $k^{\times}$,
the cyclic group $C_{p}:=\langle g\rangle$ operates trivially on
$ky$. As a result, the subalgebra $$H(C_{p})\langle y\rangle\cong
k[x,y]/(x^{p},y^{2})$$ and thus $\mathcal{C}_{H(C_{p})\langle
y\rangle}(k)=2$ by Example 2.6 (2). By
$\mathcal{C}_{H(\mathscr{G})}(k) \leq 1$ and $H(\mathscr{G})$ is
projective over $H(C_{p})\langle y\rangle$,
$\mathcal{C}_{H(C_{p})\langle y\rangle}(k)\leq 1$. It is a
contradiction.

Next, assume that $[y,y]\neq 0$. Also, if
$H(\mathscr{\underline{G}}_{red})$ is not semisimple, then similar
to the proof above one can find $g\in \underline{\mathscr{G}}_{red}$
of order $p$ such that the cyclic group $C_{p}=\langle g\rangle$
commutates with $y$. In the following, let $L:={\Lie}(\mathscr{G})$
and so $[L_{1},L_{1}]_{p}\subset N(L_{0})$ (by Lemma 3.1 (2) and
Proposition 3.2). Since $\dim_{k}L_{1}=\dim_{k}(ky)=1 $,
$[[L_{1},L_{1}]_{p},L_{1}]=0$. Now consider the quotient
$H(\underline{\mathscr{G}})\langle y\rangle/([L_{1},L_{1}])$, it
contains a sub super Hopf algebra generated by $g$ and $y$ (we
identify $g,y$ with their images in
$H(\underline{\mathscr{G}})\langle y\rangle/([L_{1},L_{1}])$). Note
that $y\not\in ([L_{1},L_{1}])$ by $[[L_{1},L_{1}]_{p},L_{1}]=0$. As
an algebra, this sub super Hopf algebra is isomorphic to
$k[x,y]/(x^{p},y^{2})$. So we also have
$\mathcal{C}_{k[x,y]/(x^{p},y^{2})}(k)=2$. A contradiction.
 \end{proof}

\begin{proposition} If $[y,y]=0$, then $\underline{H}$ is
semisimple.
\end{proposition}
\begin{proof} By Proposition 3.2, $\textbf{u}(L)$ has finite
representation type and thus it has the structure given in Lemma 3.1
(2). Assumption implies that $x_{0}=0$ and so $\textbf{u}(L_{0})$ is
semisimple. Thus $\mathscr{\underline{G}}^{\circ}$ dose not contain
a copy of $_{p}\mathbf{\alpha}_{k}$, the Frobenius kernel of the
additive group $\mathbf{\alpha}_{k}$. Then Chapter IV, Section 3
(3.7) in \cite{DG} implies that $\mathscr{\underline{G}}^{\circ}$ is
multiplicative.

Now, by Lemma 4.2,
$\underline{H}=H(\mathscr{\underline{G}}^{\circ})\#
H(\mathscr{\underline{G}}_{red})$ is semisimple too since, for
example,
gl.dim$\underline{H}=$gl.dim$H(\mathscr{\underline{G}}^{\circ})$ by
Theorem 1.1 in \cite{G}.
\end{proof}

Let $\mathscr{M}(\mathscr{\underline{G}}^{\circ})$ be the largest
multiplicative center of $\mathscr{\underline{G}}^{\circ}$ and
$_{p}\mathscr{W}(n)_{k}$ the infinitesimal group corresponding to
the restricted enveloping algebra $\textbf{u}(L_{n})$  of the
$n$-dimensional $p$-nilpotent abelian restricted Lie algebra
$L_{n}:=\oplus_{i=0}^{n-1}kx^{[p]^{i}}$ with $x^{[p]^{n}}=0$. It is
the Frobenius kernel of the $n$th \emph{Witt group}
$\mathscr{W}(n)_{k}$ (see Chapter V in \cite{DG}). Denote the $n$th
Frobenius kernel of the multiplicative group $\mu_{k}$ by
$_{p^{n}}\mu_{k}$.

\begin{proposition} If $[y,y]\neq 0$, then either $\mathscr{\underline{G}}^{\circ}$
is multiplicative or
$$\mathscr{\underline{G}}^{\circ}/\mathscr{M}(\mathscr{\underline{G}}^{\circ})\cong\;
_{p}\mathscr{W}(n)_{k}\rtimes\; _{p^{m}}\mu_{k}$$ for some $m,n\in
\mathbb{N}$.
\end{proposition}

 To show this conclusion, one preparation is needed.
 By Theorem 2.7 in \cite{FV}, $\mathscr{\underline{G}}^{\circ}/\mathscr{M}(\mathscr{\underline{G}}^{\circ})
  \cong \mathscr{U}\rtimes\; _{p^{m}}\mu_{k}$ with a $V$-uniserial
  normal subgroup $\mathscr{U}$. All $V$-uniserial groups
  are classified in \cite{FRV} and they are described as $_{p}\mathscr{W}(n)_{k}$,
$\mathscr{U}_{n,d}$ and  $\mathscr{U}_{n,d}^{j}$ respectively (See
Theorem 1 in \cite{FRV} for details). Due to the complexity of the
such groups, the Hopf structures of the algebras of measures of them
are not very clear. Incidentally, the author with his corporators
\cite{HGY} realized such Hopf structures can be described through
the path coalgebra over a loop. In fact, the coordinate rings of
such groups are denoted as $L(n,d)$ in \cite{HGY}. By definition,
for any $0\leq d\leq n$,  $L(n,d)$ is defined to  be the Hopf
algebra over $k\circlearrowleft_{p^{n}}^{c}$ (see Subsection 2.1 for
the notions) with relations:
  \begin{equation} \alpha_{p^{i}}\alpha_{p^{j}}=\alpha_{p^{j}}\alpha_{p^{i}},\;\;\textrm{for}\;0\leq i,j\leq
  n-1;\end{equation}
  \begin{equation} \alpha_{p^{i}}^{p}=0,\;\;\textrm{for}\;i< d;
  \end{equation}
   \begin{equation} \alpha_{p^{i}}^{p}=\alpha_{p^{i-d}},\;\;\textrm{for}\;i\geq
   d.
  \end{equation}

\begin{lemma} \emph{(1)} Any one of
$H(_{p}\mathscr{W}(n)_{k}),\;H(\mathscr{U}_{n,d})$ and
$H(\mathscr{U}_{n,d}^{j})$ is isomorphic to $(L(n',d'))^{\ast}$ for
some $n',d'$. In particular, $H(_{p}\mathscr{W}(n)_{k})\cong
(L(n,n))^{\ast}$.

\emph{(2)} As an algebra, there is a canonical isomorphism
$(L(n,d))^{\ast}\cong k[x]/(x^{p^{n}})$, and under such isomorphism
$x^{p^{n-d}}, x^{p^{n-d+1}},\ldots,x^{p^{n-1}}$ is a basis of the
space of primitive elements of  $L(n,d))^{\ast}$.
\end{lemma}

\begin{proof} (1) is indeed the direct consequence of the proof of
Theorem 5.1 in \cite{HGY}.

Denote the dual basis of $k\circlearrowleft_{p^{n}}^{c}$ by
$\{\alpha_{i}^{\ast}\}_{0\leq i<p^{n}-1}$. That is,

$$\alpha_{i}^{\ast}(\alpha_{j})=\left \{
\begin{array}{ll} 1 & \;\;\;\;i=j\\0 &
\;\;\;\;i\neq j.
\end{array}\right. $$
Define a map $((L(n,d))^{\ast})^{\ast}\to k[x]/(x^{p^{n}})$ through
$\alpha_{i}^{\ast} \mapsto x^{i}$ for $0\leq i<p^{n}-1$, and it is
straightforward to show this is an isomorphism of algebras. Consider
this isomorphism as an identity for short. By the relations defined
through (4.3) and (4.4), one can see that
$$\Delta(x)=1\otimes x+x\otimes 1+\sum_{i=1}^{p-1}(x^{p^{d}})^{i}\otimes (x^{p^{d}})^{p-i}+\textrm{higher items}.$$
Here ``higher items" are items such as $x^{j}\otimes x^{l}$ with
$j+l> p^{d+1}$. Therefore,
$$\Delta(x^{p^{n-d}})=(\Delta(x))^{p^{n-d}}=1\otimes x^{p^{n-d}}+x^{p^{n-d}}\otimes 1,$$
and so $\{x^{p^{n-d}}, x^{p^{n-d+1}},\ldots,x^{p^{n-1}}\}\subset
P((L(n,d))^{\ast})$. To attack that it is indeed a basis, it is
enough to show $L(n,d)$ is indeed generated by $d$ elements. In
fact, if we write $n=md+i$ for $0\leq i< d$, then relation (4.4)
shows us that
$$\alpha_{p^{(m-1)d+(i+1)}},\;\alpha_{p^{(m-1)d+(i+2)}},\;\ldots,\;\alpha_{p^{(md+i)}}$$
can generate the whole $L(n,d)$.
\end{proof}

\noindent{\textbf{Proof of Proposition 4.4.}}  Clearly, to prove the
result, there is no harm to assume that
$\mathscr{\underline{G}}=\mathscr{\underline{G}}^{\circ}$. Moreover,
one even can assume that $x_{0}\neq 0$ since otherwise
$H(\mathscr{\underline{G}})$ will be semisimple by the proof of
Proposition 4.3. Consider the quotient
$$\underline{H}(\mathscr{G})\langle y\rangle/(H(\mathscr{M}(\mathscr{\underline{G}}))^{+}).$$
We claim that $y\not\in
(H(\mathscr{M}(\mathscr{\underline{G}}))^{+})$. Otherwise,
$2y^{2}=[y,y]\in (H(\mathscr{M}(\mathscr{\underline{G}}))^{+})\cap
\textbf{u}(L_{0})$. By $H(\mathscr{G})$ is faithfully flat over
$H(\mathscr{\underline{G}})$,
$(H(\mathscr{M}(\mathscr{\underline{G}}))^{+})\cap
\textbf{u}(L_{0})=H(\mathscr{M}(\mathscr{\underline{G}}))^{+}H(\mathscr{\underline{G}})\cap
\textbf{u}(L_{0})$ which is contained in $\textbf{u}(T(L_{0})+kt)$
(see notions in Lemma 3.1 (2)) by the definition of
$\mathscr{M}(\mathscr{\underline{G}})$. By $x_{0}\in k[y,y]_{p}$,
$x_{0}\in \textbf{u}(T(L_{0})+kt)$ which is impossible. Thus
$$
\underline{H}(\mathscr{G})\langle
y\rangle/(H(\mathscr{M}(\mathscr{\underline{G}}))^{+})\cong
(\underline{H}(\mathscr{G})/H(\mathscr{M}(\mathscr{\underline{G}}))^{+}\underline{H}(\mathscr{G}))\langle
y\rangle$$ which, by Theorem 2.7 in \cite{FV}, is isomorphic to $
(H(\mathscr{U}\rtimes\; _{p^{m}}\mu_{k}))\langle y\rangle$ for some
$V$-uniserial group $\mathscr{U}$. By Lemma 4.5 (1),
$$H(\mathscr{U}\rtimes\; _{p^{m}}\mu_{k})\cong H(\mathscr{U})\# (k\mathbb{Z}_{p^{m}})^{\ast}\cong
(L(n,d))^{\ast}\# (k\mathbb{Z}_{p^{m}})^{\ast}$$ for some $n,d$ with
$d\leq n$. Owing to Lemma 4.5 (2), $\{x^{p^{n-d}},
x^{p^{n-d+1}},\ldots,x^{p^{n-1}}\}$ is  a basis of the space of
primitive elements of  $L(n,d))^{\ast}$. Therefore,
$P(L(n,d))^{\ast})=(kx^{p^{n-d}})_{p}$. Denote the Lie algebra of
$(L(n,d))^{\ast}\# (k\mathbb{Z}_{p^{m}})^{\ast}$ by $L_{0}$, then
$$L_{0}=(kx^{p^{n-d}})_{p}+kt$$
with a toral element $t$ which does not commute with $x^{p^{n-d}}$.
Now $N(L_{0})=(kx^{p^{n-d}})_{p}$. So, by Proposition 3.2 and Lemma
3.1 (2), $$[ky,ky]_{p}=(kx^{p^{n-d}})_{p}.$$ Also, by the proof of
Proposition 3.2 (part (f)), $[x,y]=0$. So, as an algebra,
$$(\ast)\;\;\;\;H(\mathscr{U})\langle y\rangle\cong k[x,y]/(x^{p^{n}},y^{2}-x^{p^{n-d}}).$$
Forming the quotient super Hopf algebra $H(\mathscr{U}\rtimes\;
_{p^{m}}\mu_{k})\langle y\rangle/(x^{p^{n-d}})$, it contains
$H(\mathscr{U})\langle y\rangle/(x^{p^{n-d}})$ as a sub super Hopf
algebra. By $H(\mathscr{U}\rtimes\; _{p^{m}}\mu_{k})\langle
y\rangle/(x^{p^{n-d}})$ has finite representation type,
$\mathcal{C}_{H(\mathscr{U})\langle y\rangle/(x^{p^{n-d}})}(k)\leq
1$. Thus $(\ast)$ implies that
$\mathcal{C}_{k[x,y]/(x^{p^{n-d}},y^{2})}(k)\leq 1$. Owing to
Example 2.6 (2), this is possible only in the case $n=d$. Thus,
$H(\mathscr{U})\cong (L(n,n))^{\ast}$ and by Lemma 4.5 (1),
$\mathscr{U}\cong\; _{p}\mathscr{W}(n)_{k}$ as desire.
$\;\;\;\;\square$

\section{Representation-finite super groups of dimension zero}

Combining the conclusions gotten in Sections 3,4, we will determine
the structure of representation finite super groups of dimension
zero in this section. The following conclusion is the direct
consequence of the proof of Proposition 2.2 (1) in \cite{LZ}.
\begin{lemma} Let $H$ be a semisimple Hopf algebra and $A$ is a
finite-dimensional twisted $H$-module algebra such that
$A\#_{\sigma}H$ exists. Then $A\#_{\sigma}H$ is of finite
representation type if $A$ is so.
\end{lemma}
The next result, which given as the Theorem 3.3 in \cite{FV2}, is
also needed.
\begin{lemma} Let $\mathscr{H}$ be an infinitesimal group such that
$\mathscr{H}/\mathscr{M}(\mathscr{H}) \cong\;
_{p}\mathscr{W}(n)_{k}$. Then $\mathscr{H}$ is commutative and
$\mathscr{H}\cong \; _{p}\mathscr{W}(n)_{k}\times
\mathscr{M}(\mathscr{H})$.
\end{lemma}

\begin{theorem} Let $\mathscr{G}$ be a finite super algebraic
$k$-group with $\mathscr{G}\neq \underline{\mathscr{G}}$ and
$H(\mathscr{G})$ be its algebra of measures. Then the following are
equivalent:

\emph{(1)} $H(\mathscr{G})$ has finite representation type.

\emph{(2)} $\emph{\textbf{u}}(\emph{Lie}(\mathscr{G}))$ has finite
representation type and either $H(\underline{\mathscr{G}})$ is
semisimple or
$$\mathscr{\underline{G}}^{\circ}/\mathscr{M}(\mathscr{\underline{G}}^{\circ})\cong\;
_{p}\mathscr{W}(n)_{k}\rtimes\; _{p^{m}}\mu_{k}$$ for some $m,n\in
\mathbb{N}$.
\end{theorem}
\begin{proof} ``(1) $\Rightarrow$ (2)" By Proposition 3.2, $\textbf{u}({\Lie}(\mathscr{G}))$ has finite
representation type. Thus there is $0\neq y\in V_{H(\mathscr{G})}$
such that $H(\mathscr{G})=H(\underline{\mathscr{G}})\langle
y\rangle$. If $[y,y]=0$, $H(\underline{\mathscr{G}})$ is semisimple
by Proposition 4.3. Otherwise, $[y,y]\neq 0$. In this case, if
$\mathscr{\underline{G}}^{\circ}/\mathscr{M}(\mathscr{\underline{G}}^{\circ})\not\cong\;
_{p}\mathscr{W}(n)_{k}\rtimes\; _{p^{m}}\mu_{k}$, then Proposition
4.4 implies that $\underline{\mathscr{G}}^{\circ}$ is
multiplicative. So together with an application of Lemma 4.2,
$H(\underline{\mathscr{G}})$ is semisimple.

``(2) $\Rightarrow$ (1)"  At first, assume that
$H(\underline{\mathscr{G}})$ is semisimple. Since $H({\mathscr{G}})$
is a super Hopf algebra, it is a Hopf algebra in the category
$^{\mathbb{Z}_{2}}_{\mathbb{Z}_{2}}\mathscr{YD}$. So
$H({\mathscr{G}})\rtimes k\mathbb{Z}_{2}$ is a usual Hopf algebra.
Lemma 2.4 (3) implies that $\textbf{u}({\Lie}(\mathscr{G}))\rtimes
k\mathbb{Z}_{2}$ is a normal sub Hopf algebra  and so we have a Hopf
surjection
$$H({\mathscr{G}})\rtimes k\mathbb{Z}_{2}\twoheadrightarrow (H({\mathscr{G}})\rtimes k\mathbb{Z}_{2})/
(\textbf{u}({\Lie}(\mathscr{G}))\rtimes
k\mathbb{Z}_{2})^{+}(H({\mathscr{G}})\rtimes k\mathbb{Z}_{2}).$$
Owing to Theorem 8.4.6 in \cite{Sch}, $H({\mathscr{G}})\rtimes
k\mathbb{Z}_{2}\cong (\textbf{u}({\Lie}(\mathscr{G}))\rtimes
k\mathbb{Z}_{2})\#_{\sigma} ((H({\mathscr{G}})\rtimes
k\mathbb{Z}_{2})/ (\textbf{u}({\Lie}(\mathscr{G}))\rtimes
k\mathbb{Z}_{2})^{+}(H({\mathscr{G}})\rtimes k\mathbb{Z}_{2}))$. By
assumption, $(H({\mathscr{G}})\rtimes k\mathbb{Z}_{2})/\\
(\textbf{u}({\Lie}(\mathscr{G}))\rtimes
k\mathbb{Z}_{2})^{+}(H({\mathscr{G}})\rtimes k\mathbb{Z}_{2})$ is
semisimple.  Note that $\textbf{u}({\Lie}(\mathscr{G}))\rtimes
k\mathbb{Z}_{2}$ has finite representation type (by Lemma 2.1) and
by Lemma 5.1, $H({\mathscr{G}})\rtimes k\mathbb{Z}_{2}$ and so
$H({\mathscr{G}})$ (by using Lemma 2.1 again) has finite
representation type.

Next, assume
$\mathscr{\underline{G}}^{\circ}/\mathscr{M}(\mathscr{\underline{G}}^{\circ})\cong\;
_{p}\mathscr{W}(n)_{k}\rtimes\; _{p^{m}}\mu_{k}$. Let
$\mathscr{N}(\underline{\mathscr{G}})$ be the nilpotent radical of
$\underline{\mathscr{G}}$. So assumption and Lemma 5.2 show that
both
$$H(\mathscr{N}(\underline{\mathscr{G}}))/\textbf{u}(L_{0})^{+}H(\mathscr{N}(\underline{\mathscr{G}}))\;\;
\textrm{and}
\;\;H(\underline{\mathscr{G}})/H(\mathscr{N}(\underline{\mathscr{G}}))^{+}H(\underline{\mathscr{G}})$$
are semisimple. By Lemma 2.5, both
$$H(\mathscr{G})/(H(\mathscr{N}(\underline{\mathscr{G}}))\langle
y\rangle)^{+}H(\mathscr{G})\;\;\textrm{and}\;\;H(\mathscr{N}(\underline{\mathscr{G}}))\langle
y\rangle/\textbf{u}(L)^{+}H(\mathscr{N}(\underline{\mathscr{G}}))\langle
y\rangle$$ are semisimple. By $\textbf{u}(L)$ is of finite
representation type and applying the same methods used in the above
paragraph twice, $H(\mathscr{G})$ has finite representation type.
\end{proof}

\section{The Auslander-Reiten quiver}

Recall an algebra $A$ is a Nakayama algebra if each indecomposable
$A$-module is uniserial. According to Theorem 2.1 in Chapter VI of
\cite{ARS}, every Nakayama algebra has finite representation type.
The main result of this section is to show that the converse is also
true for super cocommutative Hopf algebras and the Auslander-Reiten
quivers of representation-finite super cocommutative Hopf algebras
can be deduced by this result right now.

\begin{theorem}  Let $\mathscr{G}$ be a finite super algebraic
$k$-group with $\mathscr{G}\neq \underline{\mathscr{G}}$ and
$H(\mathscr{G})$ be its algebra of measures. If $H(\mathscr{G})$ is
of finite representation type, then it is a Nakayama algebra.
\end{theorem}

To show it, we begin with some observations. By the proof of
Proposition 4.1, $H(\mathscr{G})=H(\underline{\mathscr{G}})\langle
y\rangle$ which is isomorphic to
$H(\underline{\mathscr{G}}^{\circ})\#
H(\underline{\mathscr{G}}_{red})\langle y\rangle\cong
H(\underline{\mathscr{G}}^{\circ})\langle y\rangle\#
H(\underline{\mathscr{G}}_{red})$. Owing to Lemma 4.2,
$H(\underline{\mathscr{G}}_{red})$ is always semisimple. Thus Lemma
2.3 implies $H(\underline{\mathscr{G}}^{\circ})\langle y\rangle\#
H(\underline{\mathscr{G}}_{red})$ is a Nakayama algebra if and only
if $H(\underline{\mathscr{G}}^{\circ})\langle y\rangle$ is so.
Therefore, to show the theorem one can assume that
$$\underline{\mathscr{G}}=\underline{\mathscr{G}}^{\circ}.$$
Under such assumption, we have
\begin{lemma} If $H(\underline{\mathscr{G}})$ is semisimple, then
$H({\mathscr{G}})$ is a Nakayama algebra.
\end{lemma}
\begin{proof} By Nagata's Theorem (Chapter IV, $\S $ 3, 3.6),
$H(\underline{\mathscr{G}})$ is commutative. Thus
$H(\underline{\mathscr{G}})$ decomposes into a direct sum
$$H(\underline{\mathscr{G}}) =\bigoplus_{\gamma} k_{\gamma}$$ of
one-dimensional modules. Hence, we obtain
$$H({\mathscr{G}})\cong \bigoplus_{\gamma}  H({\mathscr{G}})\otimes_{H(\underline{\mathscr{G}})} k_{\gamma},$$
a direct sum of projective $H({\mathscr{G}})$-modules. Consequently,
the dimension of each projective indecomposable
$H({\mathscr{G}})$-module is bounded by 2, forcing all these modules
to be uniserial. Note that $H({\mathscr{G}})$ is a Frobenius
algebra, all projective modules are injective and vice versa. As a
result, $H({\mathscr{G}})$ is a Nakayama algebra.
\end{proof}

In the following, we always assume that
$\underline{\mathscr{G}}=\underline{\mathscr{G}}^{\circ}$ unless
stated otherwise. Using Theorem 5.3 and above lemma, we only need to
consider the case
$\underline{\mathscr{G}}/\mathscr{M}(\mathscr{\underline{G}})\cong\;
_{p}\mathscr{W}(n)_{k}\rtimes\; _{p^{m}}\mu_{k}$ for some $m,n\in
\mathbb{N}$.

\begin{lemma} If $\underline{\mathscr{G}}/\mathscr{M}(\mathscr{\underline{G}})\cong\;
_{p}\mathscr{W}(n)_{k}\rtimes\; _{p^{m}}\mu_{k}$ for some $m,n\in
\mathbb{N}$, then $H(\mathscr{M}(\mathscr{\underline{G}}))$ commutes
with $y$.
\end{lemma}
\begin{proof} If not, there exists an element $h\in
H(\mathscr{M}(\mathscr{\underline{G}}))$ such that $hy\neq yh$. By
the proof of Proposition 4.1, there is a character $\chi:
H(\underline{\mathscr{G}})\to k$ such that $y\cdot h=\chi(h)y$. Thus
assumption implies that $\chi(h)\neq \varepsilon(h)$. So$$0\neq
(\chi(h)-\varepsilon(h))y=y\cdot (h-\varepsilon(h)1)\in
(H(\mathscr{M}(\mathscr{\underline{G}}))^{+})$$ the ideal generated
by $H(\mathscr{M}(\mathscr{\underline{G}}))^{+}$. Thus $y\in
(H(\mathscr{M}(\mathscr{\underline{G}}))^{+})$ which is impossible
by the proof of Proposition 4.4.
\end{proof}

Denote by $\mathscr{B}_{0}(H(\mathscr{G}))$ the block of
$H(\mathscr{G})$ containing the trivial module $k$.

\begin{lemma}If $\underline{\mathscr{G}}/\mathscr{M}(\mathscr{\underline{G}})\cong\;
_{p}\mathscr{W}(n)_{k}\rtimes\; _{p^{m}}\mu_{k}$ for some $m,n\in
\mathbb{N}$, then

\emph{(1)} As an algebra,
$H(\mathscr{G})/(H(\mathscr{M}(\mathscr{\underline{G}}))^{+})\cong
k[y]/(y^{2p^{n}})\# (k\mathbb{Z}_{p^{m}})^{\ast}$ which is a
Nakayama algebra.

\emph{(2)} The canonical projection $\pi:\;H(\mathscr{G})\to
H(\mathscr{G})/(H(\mathscr{M}(\mathscr{\underline{G}}))^{+})$
induces an isomorphism $\mathscr{B}_{0}(H(\mathscr{G}))\cong
H(\mathscr{G})/(H(\mathscr{M}(\mathscr{\underline{G}}))^{+})$.
\end{lemma}
\begin{proof} (1) By the proof of Proposition 4.4,
$$H(\underline{\mathscr{G}})/H(\mathscr{M}(\mathscr{\underline{G}}))^{+}H(\underline{\mathscr{G}})\cong
(L(n,n))^{\ast}\#(k\mathbb{Z}_{p^{m}})^{\ast} \cong
k[x]/(x^{p^{n}})\#(k\mathbb{Z}_{p^{m}})^{\ast}$$ and
$[ky,ky]_{p}=(kx)_{p}$. Thus it is harmless to assume that $y^{2}=x$
and so
$$H(\mathscr{G})/(H(\mathscr{M}(\mathscr{\underline{G}}))^{+})\cong (L(n,n))^{\ast}\langle y\rangle \#(k\mathbb{Z}_{p^{m}})^{\ast}
\cong k[y]/(y^{2p^{n}})\# (k\mathbb{Z}_{p^{m}})^{\ast}.$$ Since $ky$
is invariant under the action of $(k\mathbb{Z}_{p^{m}})^{\ast}$, the
Jacobson radical $J_{k[y]/(y^{2p^{n}})\#
(k\mathbb{Z}_{p^{m}})^{\ast}}$ equals to $(ky)\#
(k\mathbb{Z}_{p^{m}})^{\ast}$. And so $$k[y]/(y^{2p^{n}})\#
(k\mathbb{Z}_{p^{m}})^{\ast}/J_{k[y]/(y^{2p^{n}})\#
(k\mathbb{Z}_{p^{m}})^{\ast}}\cong (k\mathbb{Z}_{p^{m}})^{\ast}$$
and $$J_{k[y]/(y^{2p^{n}})\#
(k\mathbb{Z}_{p^{m}})^{\ast}}/J_{k[y]/(y^{2p^{n}})\#
(k\mathbb{Z}_{p^{m}})^{\ast}}^{2}\cong
(k\mathbb{Z}_{p^{m}})^{\ast}.$$ From this, the Gabriel's quiver of
$k[y]/(y^{2p^{n}})\# (k\mathbb{Z}_{p^{m}})^{\ast}$ is a basic cycle
with $\dim_{k}(k\mathbb{Z}_{p^{m}})^{\ast}$ vertices. Thus it is
Nakayama.

(2) According to (1),
$H(\mathscr{G})/(H(\mathscr{M}(\mathscr{\underline{G}}))^{+})$ is
connected by noting that $x$ (and so $y$) does not commute with
$(k\mathbb{Z}_{p^{m}})^{\ast}$. It follows that the restriction
$\pi:\; \mathscr{B}_{0}(H(\mathscr{G}))\to
H(\mathscr{G})/(H(\mathscr{M}(\mathscr{\underline{G}}))^{+})$ of the
canonical projection maps the primitive central idempotent of
$\mathscr{B}_{0}(H(\mathscr{G}))$ onto the identity. Consequently,
$\pi$ is surjective. Since the ideal
$(H(\mathscr{M}(\mathscr{\underline{G}}))^{+})=H(\mathscr{G})H(\mathscr{M}(\mathscr{\underline{G}}))^{+}$
(by Lemma 6.3) is indeed generated by central idempotents not
belonging to $\mathscr{B}_{0}(H(\mathscr{G}))$, the map $\pi$ is
also injective, and our assertion follows.
\end{proof}

Let $H$ be an ordinary Hopf algebra and $M,N$ two $H$-modules. One
can equip the tensor product $M\otimes N$ with an $H$-module
structure through the comultiplication $\Delta:\;H\to H\otimes H$
and make $\Hom_{k}(M,N)$ to be an $H$-module by $(h\cdot f)(m):=\sum
h_{(1)}f(S(h_{(2)})m)$ for $f\in \Hom_{k}(M,N)$ and $h\in H$ . In
case of $H$ is a super Hopf algebra, one also can do the same
constructions by using super modules. The following result is the
counter part of Corollary 2.5 (1) in \cite{FV2} in super case.

\begin{lemma}Assume that $\underline{\mathscr{G}}/\mathscr{M}(\mathscr{\underline{G}})\cong\;
_{p}\mathscr{W}(n)_{k}\rtimes\; _{p^{m}}\mu_{k}$ for some $m,n\in
\mathbb{N}$. Let $\mathscr{B}$ be a block of $H(\mathscr{G})\rtimes
k\mathbb{Z}_{2}$ and $S,T$ be two simple modules belonging to
$\mathscr{B}$. Then there exists a character $\gamma:\;
H(\mathscr{G})\rtimes k\mathbb{Z}_{2}\to k$ such that $T\cong
k_{\gamma}\otimes S$.
\end{lemma}
\begin{proof} Note that $H(\mathscr{G})\rtimes k\mathbb{Z}_{2}$ is
an ordinary Hopf algebra. Consider the $H(\mathscr{G})\rtimes
k\mathbb{Z}_{2}$-module $\Hom_{k}(S,T)$. By Lemma 6.3,
$H(\mathscr{M}(\mathscr{\underline{G}}))$ lies in the center of
$H(\mathscr{G})\rtimes k\mathbb{Z}_{2}$. By $S,T$ belonging to the
same block, $H(\mathscr{M}(\mathscr{\underline{G}}))$ operates on
$S$ and $T$ via the same character and so acts trivially on
$\Hom_{k}(S,T)$. Hence $\Hom_{k}(S,T)$ is a $H(\mathscr{G})\rtimes
k\mathbb{Z}_{2}/(H(\mathscr{M}(\mathscr{\underline{G}}))^{+})$-module,
which is a basic algebra by Lemma 6.4 (1). Therefore,
$\Hom_{k}(S,T)$ contains a 1-dimensional submodule $k_{\gamma}$,
defined by a character $\gamma$ of $H(\mathscr{G})\rtimes
k\mathbb{Z}_{2}$. Let $\psi$ be a non-zero element of $k_{\gamma}$
and consider
$$\widehat{\psi}:\;k_{\gamma}\otimes S\to T,\;\;\psi\otimes x\mapsto \psi(x).$$
Now, for $h\in H(\mathscr{G})\rtimes k\mathbb{Z}_{2}$,
\begin{eqnarray*}
\widehat{\psi}(h\cdot (\psi\otimes x))&=&\sum
\widehat{\psi}(h_{(1)}\cdot \psi\otimes h_{(2)}\cdot x)\\
&=&\sum (h_{(1)}\cdot \psi)(h_{(2)}\cdot x)=\sum
h_{(1)}\psi(S(h_{(2)})h_{(3)}\cdot x)\\
&=&h\psi(x)=h\cdot \widehat{\psi}(\psi\otimes x).
\end{eqnarray*}
Consequently, $\widehat{\psi}$ is, as a non-zero
$H(\mathscr{G})\rtimes k\mathbb{Z}_{2}$-linear map between two
simple modules, an isomorphism.
\end{proof}

We now in the position to prove Theorem 6.1 now.\\[2mm]
\textbf{Proof of Theorem 6.1.} By Lemma 6.2 and Theorem 5.3, one can
assume that
$\underline{\mathscr{G}}/\mathscr{M}(\mathscr{\underline{G}})\cong\;
_{p}\mathscr{W}(n)_{k}\rtimes\; _{p^{m}}\mu_{k}$ for some $m,n\in
\mathbb{N}$. Owing to Lemma 6.4,
$\mathscr{B}_{0}(H(\mathscr{G}))\cong
H(\mathscr{G})/(H(\mathscr{M}(\mathscr{\underline{G}}))^{+})$ is a
Nakayama algebra.

By the remarks after Lemma 2.3, there is no harm to consider the
super modules. With such basic observation, we continue by showing
that any block $\mathscr{B}\subset H(\mathscr{G})$ containing a
1-dimensional super module is a Nakayama algebra. According to Lemma
6.5 and Lemma 2.2, every simple super $\mathscr{B}$-module is of the
form $k_{\lambda}$ for some character $\lambda:\;H(\mathscr{G})\to
k$. Given two simple super $\mathscr{B}$-modules $k_{\mu},k_{\nu}$,
we have
$$\textrm{Ext}^{1}_{H(\mathscr{G})}(k_{\mu},k_{\nu})\cong \textrm{Ext}^{1}_{H(\mathscr{G})}(k\otimes k_{\mu},k_{\nu})
\cong
\textrm{Ext}^{1}_{H(\mathscr{G})}(k,\Hom_{k}(k_{\mu},k_{\nu}))$$
which is isomorphic to
$$\textrm{Ext}^{1}_{H(\mathscr{G})}(k,k_{\nu\ast \mu^{-1}}).$$
Here the action $\ast$ denotes the convolution product and
$\mu^{-1}=\mu\ast S$. Just like in the proof of Lemma 6.5,
$H(\mathscr{M}(\mathscr{\underline{G}}))$ operates trivially on
$k_{\nu\ast \mu^{-1}}$ and so $k_{\nu\ast \mu^{-1}}$ is a module
belonging to $\mathscr{B}_{0}(H(\mathscr{G}))$ by Lemma 6.4 (2).
Thus $$\textrm{Ext}^{1}_{H(\mathscr{G})}(k,k_{\nu\ast
\mu^{-1}})\cong
\textrm{Ext}^{1}_{\mathscr{B}_{0}(H(\mathscr{G}))}(k,k_{\nu\ast
\mu^{-1}})$$ and so
$$\sum_{\nu}\dim_{k}\textrm{Ext}^{1}_{H(\mathscr{G})}(k_{\mu},k_{\nu})=
\sum_{\nu}\dim_{k}\textrm{Ext}^{1}_{\mathscr{B}_{0}(H(\mathscr{G}))}(k,k_{\nu\ast
\mu^{-1}})\leq 1$$ by $\mathscr{B}_{0}(H(\mathscr{G}))$ is a
Nakayama algebra and Theorem 9 in \cite{Hup}. Using Theorem 9 in
\cite{Hup} again, $\mathscr{B}$ is a Nakayama algebra.

Now, let $\mathscr{B}(S)$ be an arbitrary block, belonging to the
super simple module $S$. Let $T$ be another super simple
$\mathscr{B}(S)$-module and Lemma 6.5 implies that $T\cong
k_{\gamma}\otimes S$ for some character $\gamma:\;H(\mathscr{G})\to
k$. We have shown that the block corresponding to $k_{\gamma}$ is a
Nakayama algebra. By applying Theorem 2.10 in Chapter IV of
\cite{ARS}, $\Omega^{2}(k_{\gamma})=k_{\mu}$ for a suitable algebra
homomorphism $\mu:\;H(\mathscr{G})\to k$. General principles (see
Corollary 3.1.6 in \cite{Ben}) now provide isomorphisms
$$\Omega^{2}(k_{\gamma}\otimes S)\oplus \textrm{(proj)}\cong \Omega^{2}(k_{\gamma})\otimes S\cong k_{\mu}
\otimes S.$$ Since $k_{\mu} \otimes S$ is simple and not projective,
we obtain $$\Omega^{2}(T)\cong k_{\mu} \otimes S$$ is simple. By
using Theorem 2.10 in Chapter IV of \cite{ARS} again,
$\mathscr{B}(S)$ is Nakayama.$\;\;\;\;\square$

Let $\mathbf{A}_{l}$ be the type $A$ quiver of length $l$. For more
information on quivers and the definition of the stable
Auslander-Reiten quiver $\Gamma_{s}(\Lambda)$ of a self-injective
algebra $\Lambda$ the reader may consult \cite{ARS}.

\begin{corollary} Let $\mathscr{G}$ be a finite algebraic super
$k$-group, $H(\mathscr{G})$ its algebra of measures and
$\mathscr{B}\subset H(\mathscr{G})$ a block. If $H(\mathscr{G})$ has
finite representation type, then
$$\Gamma_{s}(\mathscr{B})\cong \mathbb{Z}\mathbf{A}_{(l-1)}/(\tau^{n})$$
for $l$ the Loewy length of $\mathscr{B}$ and $n$ the number of
simple $H(\mathscr{G})$-modules belonging to $\mathscr{B}$.
\end{corollary}
\begin{proof} Direct consequence of Theorem 6.1 and general result
stated in page 253 of \cite{ARS}.
\end{proof}

We use the following example to explain the results we gotten.
\begin{example} \emph{Assume that $p=3$ (although now $p<5$, it is not essential
for this example). Let $L=L_{0}\oplus L_{1}$ be a Lie superalgebra
with $L_{0}=kx+kt+kt_{1}$ and $L_{1}=ky$ with relations
$$[t,x]=[t,y]=[t,t_{1}]=0,\;\;[t_{1},x]=x+t,$$$$\;\;[t_{1},y]=2y,\;\;[x,y]=0,\;\;[y,y]=x+t.$$}
\emph{The p-mapping is given by
$$t^{[p]}=t,\;\;t_{1}^{[p]}=t_{1},\;\; x^{[p]}=0.$$ By Lemma 3.1
(2), $\textbf{u}(L)$ has finite representation type. Let
$e_{0}:=1-t^{2},\;e_{1}:=2t+2t^{2},\;e_{2}:=t+2t^{2}$, then
$$\textbf{u}(L)=\textbf{u}(L)e_{0}\oplus\textbf{u}(L)e_{1}\oplus
\textbf{u}(L)e_{2}$$ is the block decomposition of $\textbf{u}(L)$.}

\emph{For $\textbf{u}(L)e_{0}$, by $te_{0}=t-t^{3}=0$, it is
isomorphic to
$k\{y,t_{1}\}/(y^{6},t_{1}^{3}-t_{1},t_{1}y-yt_{1}-2y)$. Note the
the subalgebra generated by $t_{1}$ is isomorphic to
$(k\mathbb{Z}_{3})^{\ast}$ and the subalgebra generated by $y$ is a
$(k\mathbb{Z}_{3})^{\ast}$-module algebra through the action
$t_{1}\cdot y:=[t_{1},y]=2y$. Thus $\textbf{u}(L)e_{0}\cong
k[y]/(y^{6})\# (k\mathbb{Z}_{3})^{\ast}$. It is not hard to see that
the group algebra of the largest multiplicative center is the
algebra generated by $t$ and thus we indeed have
$\textbf{u}(L)/(t)\cong \textbf{u}(L)e_{0}$. All facts stated in
Lemma 6.4 are verified in this case.}

\emph{For $\textbf{u}(L)e_{1}$, by $te_{1}=e_{1}$, it is isomorphic
to
$k\{x,y,t_{1}\}/(y^{2}-x-1,x^{3},t_{1}^{3}-t_{1},t_{1}y-yt_{1}-2y)$.
We will show that it is a super simple algebra. This is equivalent
to show that $\textbf{u}(L)e_{1}\# k\mathbb{Z}_{2}$ is simple.
Indeed, denote the generator of $\mathbb{Z}_{2}$ by $g$ and define
$$I_{2,-1}:=\left (
\begin{array}{cc} 1&0\\0&-1
\end{array}\right)\;\;\textrm{and}\;\;I_{1,2}:=\left (
\begin{array}{cc} 0&1\\1&0
\end{array}\right).$$
As usual, let $I_{2}$ be the $2\times 2$ identity matrix. Consider
the map
$$\phi:\;\textbf{u}(L)e_{1}\# k\mathbb{Z}_{2}\to \mathbf{M}_{6}(k)$$
by sending
$$t_{1}\mapsto \left (
\begin{array}{ccc} 0&0&0
\\                 0&I_{2}&0
\\                 0&0&2I_{2}
\end{array}\right),\;\;y\mapsto \left (
\begin{array}{ccc} 0&I_{1,2}&0
\\                 0&0&I_{1,2}
\\                 I_{1,2}&0&0
\end{array}\right),$$
$$x\mapsto \left (
\begin{array}{ccc} -I_{2}&0&I_{2}
\\                 I_{2}&-I_{2}&0
\\                 0&I_{2}&-I_{2}
\end{array}\right),\;\;g\mapsto \left (
\begin{array}{ccc} I_{2,-1}&0&0
\\                 0&I_{2,-1}&0
\\                 0&0&I_{2,-1}
\end{array}\right)$$
By direct computations, one can show $\phi$ is an algebra
isomorphism.}

\emph{Similarly, $\textbf{u}(L)e_{2}\# k\mathbb{Z}_{2}\cong
\mathbf{M}_{6}(k)$ too. Thus $\textbf{u}(L)$ is a Nakayama algebra.}

\end{example}
We end this section with the following remarks.
\begin{remark} It is known that for an infinitesimal group
$\underline{\mathscr{G}^{\circ}}$,
$H(\underline{\mathscr{G}}^{\circ})$ has finite representation type
if and only if it is a Nakayama algebra (see Theorem 2.7 \cite{FV}).
But for a constant group $G$, $kG$ may be not a Nakayama algebra
even $kG$ is of finite representation type. That is, for a
representation-finite finite algebraic $k$-group
$\underline{\mathscr{G}}$, $H(\underline{\mathscr{G}})$ may not be a
Nakayama algebra. Contrast to the ordinary finite algebraic groups,
Theorem 6.1 tells us that the phenomenon will not appear in super
case.

\end{remark}

\section*{Acknowledgements} The author is supported by
Japan Society for the Promotion of Science under the item ``JSPS
Postdoctoral Fellowship for Foreign Researchers" and Grant-in-Aid
for Foreign JSPS Fellow. I would gratefully acknowledge JSPS. I
would like thank Professor A. Masuoka for stimulating discussions
and his encouragements.

\end{document}